\documentclass{article}
\usepackage{amssymb}
\usepackage{amsmath}

\let\<\langle
\let\>\rangle
\font\bigcmsy=cmsy10.pk scaled 2000
\def\bigtimes{\mathop{\,\vrule width0pt depth2pt height8pt
            \smash{\lower2pt\hbox{\bigcmsy\char'002}}\,}\limits}

\begin{document}

\begin{center}
\Large{\textbf{Order separability of HNN-extensions and free products with commutative subgroups.}}
\end{center}

\begin{center}
\textbf{Vladimir V. Yedynak}
\end{center}

\begin{abstract}
This paper is devoted to the investigation of the property of order separability for HNN extensions and free products with commutative subgroups. Particularly it was proven that HNN extension of a free group with maximal connected cyclic subgroups is 2-order separable.


\end{abstract}

\section{Introduction.}

Definition. Consider a group $G$ and a  natural number $n$ which is greater than 1. We say that $G$ is $n$-order separable if for each elements $g_1,..., g_n$ such that $g_i$ is not conjugate to $g_j^{\pm1}$ for each $i\neq j$ there exists a homomorphism $\varphi$ of $G$ onto a finite group so that the orders of $\varphi(g_1),..., \varphi(g_n)$ are pairwise different.

The property of 2-order separability was investigated in [1] for free groups. Later it was generalized in [2] by Daniel T. Wise who proved that free groups are omnipotent.

Definition. A group $G$ is said to be omnipotent if for each elements $g_1,..., g_n$ of $G$ such that no two of them have nontrivial conjugate powers there exists a number $k$ such that for any ordered sequence of natural     numbers $l_1,.., l_n$ there exists a homomorphism $\varphi$ of $G$ onto a finite group such that the order of $\varphi(g_i)$ equals $kl_i$.

In [5] in was shown that 2-order separability is inherited by free products.

Notice that the property of 2-order separability for free products of groups was used in [3]. The property of 2-order separability strengthens the residual finiteness. Therein after we shall investigate HNN extensions. The subgroup separability of HNN extensions was investigated by Graham A. Niblo in [4]

The main results of this paper are the following.

\textbf{Theorem 1}. HNN extension of a free group with maximal connected cyclic subgroups is 2-order separable.

\textbf{Theorem 2}. Consider groups $A$ and $B$ which are subgroup separable and 2-order separable. Consider the group $G$ which is the free product of $A$ and $B$ with commutative subgroups $M$ and $N$ that is $G$ is the quotient group of $A\ast B$ on the normal closure of the mutual centralizer of $M$ and $N$ in $A\ast B$. Besides $A$ is $M$-separable and $B$ is $N$-separable that is $M$ is the intersection of finite index subgroups of $A$ and $N$ is the intersection of finite index subgroups of $B$. Consider elements $u$ and $v$ of $G$ such that $u$ is not conjugate to $v^{\pm1}$ and the following conditions is not held: $u=x^{-1}mxy^{-1}ny, v=x'^{-1}mx'y'^{-1}ny', x, x'\in A, y, y'\in B, m\in M, n\in N$. Then there exists a homomorphism of $G$ onto a finite group such that the images of $u$ and $v$ have different orders.

Definition. A group $G$ is called subgroup separable if each finitely generated subgroup of $G$ coincides with the intersection of finite index subgroups of $G$.

\section{Auxiliary notations.}

Consider the group $G$ generated by the set $\{g_i| i\in I\}$ and acting on the set $X$. Construct the graph $\Gamma$ such that the set of vertices of $\Gamma$ coincides with $X$. If $p\circ g_i=q$ and $p\circ g_i^{-1}=r$ for some $p, q, r\in X$ and for some $i\in I$ then the vertices $p$ and $q$ are connected by the edge going from $p$ into $q$ labelled by $g_i$ and $p$ and $r$ are connected by the edge going from $r$ into $p$ labelled by $g_i$. We consider that for each edge $f$ of $\Gamma$ with a label there exists the inverse edge $f'$ without a label such that the begin point of $f$ coincides with the end of $f'$ and the beginning of $f'$ is the end of $f$. The graph $\Gamma$ is oriented and edges with labels are positively oriented edges of $\Gamma$ and the rest edges are negatively oriented edges of $\Gamma$.

Definition 1. The graph $\Gamma$ which is constructed as above is called the action graph of the group $G$ with respect to its generating set $\{g_i| i\in I\}$.

We shall omit further the references about the generating set if it is fixed. Besides throughout the paper considering the action graphs of free product we shall always suppose that the set of generators coincides with the union of free factors. And when we shall deal with HNN extension the set of generators will be the union of elements from the base group and an element corresponding to a stable letter.

Consider the action graph $\Gamma$ of the group $G=\<g_i| i\in I\>$.

The label of the edge $e\in\Gamma$ is denoted as Lab $(e)$. The label of the path $S=e_1...e_n$ of $\Gamma$ is the element of the group $G$ equal to $\prod_{i=1}^n$ Lab $(e_i)'$ where Lab $(e_i)'$ equals Lab $(e_i)$ in case $e_i$ is positively oriented or Lab $(e_i')^{-1}$ otherwise where $e_i'$ is the edge inverse to $e_i$. The beginning and the end of the edge $e$ will be denoted as $\alpha(e)$ and $\omega(e)$ correspondingly while $\alpha(S), \omega(S)$ are the beginning and the end of a path $S$.

Definition 2. Consider the element $u$ of the group $G$ which does not equal to unit. Let $\Gamma$ be the action graph of the group $G=\<g_i|i\in I\>$. Fix the vertex $p$ in $\Gamma$. Then the \textsl{$u$-cycle $T$} of the graph $\Gamma$ going from $p$ is the set of paths $\{S_i=\{e_j^i| j\in J_i\}|i\in J\}$ satisfying the following conditions:

1) $\alpha(S_l)=p$;

2) there exists a one to one correspondence between the paths $S_l$ and the presentations of $u$ in generators $\{g_i|i\in I\}: u=g_{i_1}^{\varepsilon_1}...g_{i_k}^{\varepsilon_k}, $
min$(|i_j-i_{j+1}|, \varepsilon_j+\varepsilon_{j+1})>0,
\varepsilon_i\in\{-1; 1\}, 1\leqslant j\leqslant k, 1\leqslant
i\leqslant k$. Besides if $\varepsilon_j=1$ then edge $e_{kn+j}^l$ is positively oriented and its label equals $g_{i_j}$ for each natural $n$; if $\varepsilon_j=-1$ then $e_{kn+j}^l$ is negatively oriented and the label of the edge which is inverse to $e_{i_j}$ equals $g_{i_j}$ for each natural $n$; in case the path $S_l$ is finite and is composed of $r$ edges we consider that indices are modulo $r$;

3) there is no closed subpath $K$ in the path $S_l$ which differs from $S_l$ and satisfies the conditions 1), 2);

4) the path $S_l$ is the maximal path on the entry which satisfies the conditions 1)-3);

The paths which compose a $u$-cycle will be called the representatives of a $u$-cycle.

Definition 3. Suppose that some representative of the $u$-cycle $T$ is composed from the finite set of edges and has the label $u^k$. Then we shall say that the \textsl{length} of the $u$-cycle $T$ is equal to $k$.

Suppose we have a group $G$ which is defined by its presentation $G=<x_i, i\in I| R_j, j\in J>$. It is obvious that a graph $\Gamma$ is the action graph of the group $G=\<x_i| i\in I\>$ if and only if it satisfies the following conditions:

1) for each vertex $p$ of $\Gamma$ and for each $i\in I$ there exist exactly one edge with label $x_i$ going from $p$ and exactly one edge going into $p$ labelled by $x_i$;

2) for each $j\in J$ and for each vertex $p$ of $\Gamma$ the $R_j$-cycle of $\Gamma$ going from $p$ has length 1;

Using this remark we can describe action graphs of the HNN extension $H$ of a group $G$ with connected subgroups $A, B$ and stable letter $t$ and isomorphism $\varphi:A\rightarrow B$ which are so that $H$ acts freely on the set of vertices of a graph. The action graph $\Gamma$ of such kind satisfies the following properties.

1) for each vertex $p$ of $\Gamma$ and for each $x\in G\cup\{t\}$ there exists a singular edge with label $x$ going into $p$ and there exists exactly one edge with label $x$ going from $p$;

2) for each vertex $p$ from $\Gamma$ the maximal connected subgraph containing $p$ whose positively oriented edges are labelled by the elements of $G$ is the Cayley graph of $G$ with the generating set $\{G\}$;

3) for each edge $f$ with label $a\in A\setminus\{1\}$ there exists the edge $g$ with label $\varphi(a)$ and there exist edges $t_1, t_2$ with label $t$ such that $\alpha(t_1)=\alpha(f), \alpha(t_2)=\omega(f), \omega(t_1)=\alpha(g), \omega(t_2)=\omega(g)$;

Definition 4. Consider two vertices $p$ and $q$ of the graph $\Gamma$ which belong to one component of $\Gamma$. Then the
\textsl{distance} between $p$ and $q$ is the minimal length of the finite path which connects $p$ and $q$. The length of the finite path is the number of edges which compose this path.

Definition 5. Consider the graph $\Gamma$, the cycle $S=e_1...e_n$ in it and the nonnegative integer number $l$. Wa say that $S$ does not have \textsl{$l$-near vertices} if for each $i, j, 1\leqslant i<j\leqslant n$ the distance between the vertices
$\alpha(e_i), \alpha(e_j)$ is greater or equal than min$(|i-j|, n-|i-j|, l+1)$.

Note that when $l=0$ this condition means that $S$ does not have self-intersections. In this case we shall say that $S$ is simple. When $l=1$ we shall say that $S$ does not have near vertices. Remark also that if $S$ does not have $l$-near vertices then $S$ does not have $k$-near vertices for each nonnegative integer $k$ which is less than $l$.

We shall use the notation $A(p)$ to denote the subgraph containing the vertex $p$ which is the Cayley graph of the subgroup $A$ with generating set $\{A\}$. The subgraph $B(p)$ is defined analogically.

Let's characterize the action graphs of the group $R=A\ast B/\<\<[M, N]\>\>^{A\ast B}, M<A, N<B$ which corresponds to action of $R$ on the set in which $A$ and $B$ act freely. Such graph $\Lambda$ satisfies the following properties.

1) for each vertex $p$ of $\Lambda$ and for each $g\in A\cup B$ there exists exactly one positively oriented edge with label $g$ going from $p$ and exactly one edge with label $g$ going into $p$;

2) for each vertex $p$ from $\Lambda$ the maximal connected subgraph $A(p)$ containing $p$ whose positively oriented edges are labelled by elements of $A$ is the Cayley graph of the group $A$ with the generating set $\{A\}$; the same condition is true for the group $B$ and the subgraph $B(p)$;

3) for each vertex $p$ from $\Lambda$ the maximal connected subgraph containing $p$ whose positively oriented edges are labelled by the elements of $M\cup N$ is the Cayley graph of the group $M\times N$ with the set of generators $\{M\cup N\}$;

The maximal connected subgraph containing the vertex $p$ from $\Lambda$ whose positively oriented edges are labelled by the elements of $M$ ($N$) will be denoted as $M(p)$ ($N(p)$ correspondingly).

The defined graphs will be referred to as the free action graphs of corresponding groups.

In what follows we shall use the following notation. Consider the graph $\Gamma$ and its subset $T$. Consider $t$ copies of $\Gamma: \Gamma_1,..., \Gamma_t$. Then $T^i$ is the set of $\Gamma_i$ which corresponds to $T$ in $\Gamma$.

The homomorphism corresponding to the action of a group on the set of vertices of an action graph $\Gamma$ will denoted as $\varphi_{\Gamma}$.

\section{Proof of theorem 1.}

Put $\tilde{F}=<F, t; t^{-1}at=b>$ where $F$ is a free group, $t$ is a stable letter, $a, b$ are the elements of $F$ which are not proper powers in $F$. Consider the elements $u$ and $v$ such that $u$ is conjugate to neither $v$ nor $v^{-1}$. We may suppose that $u$ and $v$ are cyclically reduced.

Consider that $u, v\in F$. If one of the elements say $v$ does not belong to any connected subgroup then as it was stated in [6] there exists a homomorphism $\psi$ of $F$ onto a finite group such that $|\psi(a)|=|\psi(b)|=|\psi(u)|>|\psi(v)|$. According to the properties of $\psi$ there exists the homomorphism $\tilde{\psi}$ of the group $\tilde{F}$ onto the HNN extension of the group $\psi(F)$ with connected subgroups $\<\psi(a)\>, \<\psi(b)\>$ and isomorphism $\varphi': \<\psi(a)\>\rightarrow\<\psi(b)\>, \varphi': \psi(a)\mapsto\psi(b)$ such that $\tilde{\psi}|_F=\psi$. Since the group $\tilde{\psi}(F)$ is virtually free and the orders of $\tilde{\psi}(a), \tilde{\psi}(b)$ are different and finite there exists the homomorphism $\psi'$ of the group $\tilde{\psi}(F)$ onto a finite group such that $|\tilde{\psi}(u)|=|\psi'(\tilde{\psi}(u))|, |\tilde{\psi}(v)|=|\psi'(\tilde{\psi}(v))|$. Suppose now that
$u, v\in\<a\>, \<b\>$. Then we may consider that $u=a^n, v=b^m, |n|\neq|m|$. In this case it is possible to consider that there exists a prime number $p$ such that for some nonnegative integer $t'$ $p^{t'}\mid n, p^{t'}\nmid m$. According to [6] there exists a homomorphism $\psi_1$ of $F$ onto a finite $p$-group such that $|\psi_1(a)|=|\psi_1(b)|>nm$. As in the previous case there exists a homomorphism $\tilde{\psi_1}$ of $\tilde{F}$ onto HNN extension $\tilde{F}'$ of the group $\psi_1(F)$ with connected subgroups $\<\psi_1(a)\>, \<\psi_1(b)\>$ and isomorphism $\varphi_1: \<\psi_1(a)\>\rightarrow\<\psi_1(b)\>, \varphi_1: \psi_1(a)\mapsto\psi_1(b)$ and stable letter $\psi_1(t)$. Due to conditions on $n, m, |\psi_1(a)|$ the orders of images of $a^n, a^m$ are different. Now it is sufficient to use the residual finiteness of the group $\tilde{F}'$.

Henceforth we put $A=\<a\>, B=\<b\>$.

Consider that $u\notin F, v\in F$. Fix some reduced notation for the element $u$: $u=g_0\prod_{i=1}^nt^{\varepsilon_i}g_i, n\geqslant1$. Fix an arbitrary prime $p$. In this notation there are no entries of the kind $t^{-1}a't, tb't^{-1}, a'\in A, b'\in B$. Since $A, B$ are cyclic and $p'$-isolated $F$ is residually $p$-finite with respect to the membership in the subgroups $\<a\>, \<b\>$. So there exists the homomorphism $\psi_1$ of $F$ onto a finite $p$-group such that $\{\psi_1(g_0), \psi_1(g_1),..., \psi_1(g_n)\})\setminus\{\psi_1(A)\cup\psi_1(B)\}=\psi_1(\{g_0, g_1,..., g_n\}\setminus\{A\cup B\})$. Put $n=$ max $(|\psi_1(b)|, |\psi_1(a)|)$. Using lemma 1 from [6] we state that there exists the homomorphism $\psi_2$ of $F$ onto a finite $p$-group such that $|\psi_2(a)|=|\psi_2(b)|>N$. Consider the homomorphism $\psi': F\rightarrow\psi_1(F)\times\psi_2(F), \psi': f\mapsto\psi_1(f)\psi_2(f)$. Notice that $\{\psi'(g_0), \psi'(g_1),..., \psi'(g_n)\})\setminus\{\psi'(A)\cup\psi'(B)\}=\psi'(\{g_0, g_1,..., g_n\}\setminus\{A\cup B\}), |\psi'(a)|=|\psi'(b)|$. So there exists the homomorphism $\psi'': \tilde{F}\rightarrow\tilde{F}'=<\psi'(F), \psi'(t)| \psi'(t)^{-1}\psi'(a)\psi'(t)=\psi'(b)>$ which is so that $\psi''(u)$ does not belong to a subgroup which is conjugate to $\psi''(F)$. Hence the order of $\psi''(u)$ is infinite. Since the group $\tilde{F}'$ is virtually free and the orders of elements $\psi''(u), \psi''(v)$ are different it is sufficient to use the residual finiteness of the group $\tilde{F}'$.

Before considering the general case we change notations: we consider elements $u'\notin F, v'\notin F$. As usual $u'$ and $v'$ are cyclically reduced. Fix the reduced notations for elements $u', v': u'=g_0\prod_{i=1}^mt^{\varepsilon_i}g_it^{\varepsilon_{m+1}}, v'=h_0\prod_{j=1}^nt^{\delta_j}h_jt^{\delta_{n+1}}, m, n\geqslant0, \varepsilon_i, \delta_j\neq0$. In these notations there are no the following entries: $t^{-1}g_it, t^{-1}h_jt, tg_kt^{-1}, th_lt^{-1}, 1\leqslant i, k\leqslant n+1, 1\leqslant j, l\leqslant m+1, g_i, h_j\in A, g_k, h_l\in B$. Define the set $\Delta=\{g_i, h_j, g_ig_k^{-1}, h_jh_l^{-1}, g_ih_j^{-1}, x_s^{-1}qx_s| s\in I, 1\leqslant i, k\leqslant n+1, 1\leqslant j, l\leqslant m+1\}$ and an arbitrary prime number $p$. The element $q=1$ if $u'$ and $v'$ do not belong to one coset on some connected subgroup of $\tilde{F}$. Otherwise in case $v'=u'q, q=v'u'^{-1}$. Besides the elements $x_s, s\in I$ coincide with all possible beginnings of the word $u'$. Note that according to the definition of the irreducible notation [7] at least one of the elements $x_k^{-1}qx_k$ do not belong to $A\cup B$ if $q\neq1$.

Consider an arbitrary prime $p$. Since $A, B$ are cyclic and $p'$-isolated there exists the homomorphism $\psi_1$ of $F$ onto a finite $p$-group such that $\psi_1(\Delta\setminus\{A, B\})=\psi_1(\Delta)\setminus\psi_1(\{A, B\}), \psi_1(\Delta\setminus\{1\})=\psi_1(\Delta)\setminus\{1\}$. Put $N=$ max $(|\psi_1(a)|, |\psi_1(b)|)$. In [6] it is proven that there exists the homomorphism $\psi_2$ of $F$ onto a finite $p$-group such that $|\psi_2(a)|=|\psi_2(b)|>N$. The homomorphism $\psi_1\times\psi_2: F\rightarrow\psi_1(F)\times\psi_2(F), \psi_1\times\psi_2: f\mapsto\psi_1(f)\psi_2(f)$ maps $\tilde{F}$ onto HNN extension of the group $\psi_1\times\psi_2(F)$ with connected subgroups $\<\psi_1\times\psi_2(a)\>, \<\psi_1\times\psi_2(b)\>$. Due to the properties of $\psi_1$ and the conjugacy theorem for HNN extensions [7] the elements $\psi_1\times\psi_2(u), \psi_1\times\psi_2(v)$ do not belong to a group conjugate with $\psi_1\times\psi_2(F)$ and $\psi_1\times\psi_2(u)$ is not conjugate with $\psi_1\times\psi_2(v)^{\pm1}$. Thus we may consider that we deal with elements $u=\psi_1\times\psi_1(u'), v=\psi_1\times\psi_1(v')$ from HNN extension $\tilde{G}=<G, t| t^{-1}at=\varphi(a), a\in A>$ of the finite group and $u$ is not conjugate to $v^{\pm1}$. Besides $u, v$ do not belong to subgroups which are conjugate with the base subgroup of $\tilde{G}$. Put also $\psi_1\times\psi_2(g_i)=g_i', \psi_1\times\psi_2(h_j)=h_j',1\leqslant i\leqslant n, 1\leqslant j\leqslant m$.

Define the following notation. Let $\Gamma$ be the free action graph of some HNN extension $\tilde{H}=<H, t| t^{-1}at=\varphi(a), a\in A, \varphi(a)\in B, A\rightarrow^{\varphi}B>$. Fix in $\Gamma$ the Cayley graph $K$ of the group $B$ with the generating set $\{A\}$. For each natural number $n, n>1,$ we construct the new free action graph $\delta_n(\Gamma, K)$ of the group $\tilde{H}$ from $\Gamma$ in the following way. Consider $n$ copies of the graph $\Gamma: \Delta_1,..., \Delta_n$. For each $i, 1\leqslant i\leqslant n$ and for each edge $h^i$ with label $t$ going into the vertex from $K^i$ we delete the edge $h^i$ and connect vertices $\omega(h^i)$ and $\alpha(h^{i+1})$ by the new edge with label $t$ going from $\alpha(h^{i+1})$ (indices are modulo $n$). By the analogy we can construct the graph $\delta_n(\Gamma, R)$ considering the Cayley graph of $A$ instead of the Cayley graph of $B$. In this case we delete every edge $f_i$ with label $t$ going from the vertex of the subgraph $R^i$ and vertices $\alpha(f^{i+1})$ and $\omega(f^i)$ are connected by the edge with label $t$ going from the vertex $\alpha(f^{i+1})$. In both cases we obtain the free action graphs of HNN extensions.

Consider paths $S=e_1...e_k, T=f_1...f_l$ of a free action graph of HNN extension. In each path define the maximal set of edges $S'=\{e_{k_1},..., e_{k_n}\}, T'=\{f_{l_1},..., f_{k_m}\}$ such that for each $i, j$ the label of $e_{k_i}$ or its inverse is $t$ and the label of $f_{l_j}$ or its inverse is $t$. We shall say that $S$ and $T$ are $G$-near or that $S$ is $G$-near to $T$ if the following conditions are true:

1) sets $S'$ and $T'$ are nonempty and consist of the same number of elements ($k_1<...<k_n, l_1<...<l_m, m=n$)

2) for each $j=1,..., n, $ if the edge $e_{k_j}$ is positively oriented then the edge $f_{l_j}$ is also positively oriented; besides $A(\alpha(e_{k_j}))=A(\alpha(f_{l_j})), B(\omega(e_{k_j}))=B(\omega(f_{l_j}))$; if the edge $e_{k_j}$ is negatively oriented then the edge $f_{l_j}$ is also negatively oriented and $B(\alpha(e_{k_j}))=B(\alpha(f_{l_j})), A(\omega(e_{k_j}))=A(\omega(f_{l_j}))$

Since $\tilde{G}$ is virtually free it is cyclic subgroup separable [2]. Thereby there exists the homomorphism $\psi$ of $\tilde{G}$ onto a finite group such that the representatives $u$- and $v$-cycles in $Cay(\psi(\tilde{G}); \psi(G)\cup\{\psi(t)\})$ have no 3-near vertices. Besides due to the residual finiteness of the group $\tilde{G}$ it is possible to consider that all words from $\tilde{G}$ whose length is less or equal than 10 have nonunit image. Suppose that $|\psi(u)|=|\psi(v)|$. Consider the Cayley graph $\Gamma$ of the group $\psi(\tilde{G})$ with the generating set $\psi(\{G, t\})$. Consider the inductive process of constructions of free action graphs $\Gamma_k$ of the group $\tilde{G}, 0\leqslant k\leqslant M$. The symbol $M$ is either a natural number of the infinite symbol. It will be chosen later. Graphs $\Gamma_k$ satisfy the following properties:

1) the length of each $u$-cycle divides the length of a maximal $u$-cycle; the analogical condition is true for $v$-cycles;

2) lengths of maximal $u$- and $v$-cycles coincide;

3) the representatives of $u$ and $v$-cycles have no 3-near vertices

4) in the graph $\Gamma_k$ in case $k>0$ there exists the path $S_k$ whose length is greater or equal than $k$ such that there exists the representative of a maximal $u$-cycle containing the subpath which is $G$-near to the path $S_k$ and each representative of a maximal $v$-cycle has a subpath $G$-near to $S_k$;

5) if $f_k$ is a last edge of a path $S_k$ then the label of $f_k$ or its inverse is $t$

Notice that the absence of 3-near vertices in a cycle involves the absence of 2- and 1-near vertices in it.

Put $\Gamma_0=\Gamma$. The properties 1), 2) 3) are obviously true for the graph $\Gamma_0$ according to the definition of the Cayley graph and properties of $\psi$. Suppose there exists the graph $\Gamma_k, k\geqslant0$. Construct the graph $\Gamma_{k+1}$. Put $p=\omega(S_k)$ when $k>0$. If $k=0$ we consider that $p$ is an arbitrary vertex of $\Gamma_k$.

Put $S_k=e_1...e_{l_k}$. Consider the representative of a maximal $u$-cycle which has the subpath $G$-near to $S_k$ and corresponding the chosen notation for $u: T=f_1...f_r$. Consider also the subpath $f_{l+1}...f_{l+l_k}$ which is $G$-near to $S_k$. If $k=0$ then $T$ is an arbitrary representative of $u$-cycle passing through $p$ and corresponding to the notation of $u$ mentioned above, besides we put also $l=l_k=0$.

Consider that Lab$f_{l+l_k+1}\in G$. Due to the definition of the reduces notation for HNN extensions [7] the edge $f_{l+l_k+2}$ or its inverse has label $t$. We may also consider that the last edge of $S_k$ or its inverse has label $t$. Let $n$ be the length of a maximal $u$-cycle in $\Gamma_k$. Construct the graph $\Gamma_{k+1}'=\delta_n(\Gamma_k, R)$ where $R=A(\omega(f_{l+l_k+1}))$, in case $f_{l+l_k+2}$ is positively oriented or $R=B(\omega(f_{l+l_k+1}))$ otherwise. The graph $\Gamma_{k+1}'$ satisfies the property 3) because otherwise this property is not held for $\Gamma_k$. The property 1) is true for $\Gamma_{k+1}'$ because of the condition on the number of copies of $\Gamma_k$ and of the property 1) for $\Gamma_k$. Consider the path $S_{k+1}'=S_k^1\cup e_{l_k+1}'^1\cup e_{l_k+2}'^1$. Here the edge $e_{k_l+1}'$ is an arbitrary edge with label from $G$ going into the vertex of the subgraph $R^1$ and going away from $\omega(e_{l_k}^1)$. The edge $e_{l_k+2}'^1$ goes away from from the end of the edge $e_{k_l+1}'^1$ and goes into the vertex of the second copy of the graph $\Gamma_k$. Notice that the path $S_{k+1}'$ has length greater or equal than $k+1$ and according to the algorithm of construction of the graph $\delta_n(\Gamma_k, R)$ the edge $e_{l_k+2}'^1$ or its inverse has label $t$. According to the algorithm we may conclude that the vertex $\alpha(T)^1$ is the beginning of the maximal in $\Gamma_{k+1}'$ $u$-cycle that has the subpath $G$-near to the path $S_{k+1}'$. Let $Y$ be a representative of an arbitrary maximal $v$-cycle in $\Gamma_k$ and $Y$ corresponds to the chosen notation of $v$. Put $q=\alpha(Y), Y=g_1...g_r$. The path $g_{j+1}...g_{j+l_k}$ is $G$-near to $S_k$. Consider the representative $Y'$ of a $v$-cycle from $\Gamma_{k+1}'$ going from the vertex $q^1$ and corresponding to the chosen notation for $v$. It has the subpath $g_{j+1}^1...g_{j+l_k}^1$ which is $G$-near to the path $S_k^1$. Denote by $D^1$ the subgraph $B(\omega(e_{l_k+2}'^1))$ if $R$ is the Cayley graph of $A$ or $D^1=A(\omega(e_{l_k+2}'^1))$ in case $R$ is the Cayley graph of $B$. If $g_{j+l_k+1}$ does not go into the vertex of $R^1$ then the vertices of $Y'$ belong to the set of vertices of the first copy of the graph $\Gamma_k$. Indeed if $Y'$ has as vertex of another copy of $\Gamma_k$ then $Y'$ contains the edge $g_s^1$ going into the vertex of $R^1\cup D^1$. If $g_s^1$ into the vertex of $R^1$ then there exist at least 7 edges in $Y'$ between edges $g_s^1, g_{j+l_k+1}^1$ because all reduces words from $\tilde{F}$ whose length less or equal than 10 have nonunit images after $\psi$. Hence $Y'$ has 1-near vertices. If $g_s^1$ goes into the vertex of $D^1$ then there exist at least 6 edges in $Y'$ between edges $g_s^1, g_{j+l_k+1}^1$. And in this case $Y'$ has 2-near vertices. Due to the obtained violation and conditions 1), 2) for $\Gamma_k$ the length of the $v$-cycle corresponding to $Y'$ is less than the length of a maximal in $\Gamma_{k+1}'$ $u$-cycle. In this case we have the required homomorphism $\varphi_{\Gamma_{k+1}'}$ and $M=k$. Suppose now that $g_{j+l_k+1}^1$ goes into the vertex from $R^1$. The next edge $g_{j+l_k+2}^1$ of $Y'$ or its inverse has label $t$. If its orientation does not coincide with the orientation of $f_{l+l_k+2}^1$ then $\omega(g_{j+l_k+2}^1)$ belongs to the first copy of $\Gamma_k$. Besides all vertices of $Y'$ belong to the first copy of $\Gamma_k$ because otherwise there is an edge $g_q^1$ of $Y'$ going into $R^1\cup D^1$. If $g_q^1$ goes into the vertex of $R^1$ then $Y'$ has 1- or 0-near vertices in $R^1$ and if $g_q^1$ goes into the vertex of $D^1$ then $Y'$ has 1-near vertices in $D^1$. This is established in the same way as in the case when we considered the edge $g_{j+l_k+1}^1$. Thus in case $f_{l+l_k+2}^1$ and $g_{j+l_k+2}^1$ have different orientation the length of a maximal $u$-cycle in $\Gamma_{k+1}'$is greater than the length of an arbitrary maximal $v$-cycle in $\Gamma_{k+1}'$ and the theorem is proven letting $M=k$. So we may assume that $f_{l+l_k+2}^1$ and $g_{j+l_k+2}^1$ have coinciding orientations and hence $Y'$ possesses the subpath $G$-near to $S_{k+1}'$ and we consider that $\Gamma_{k+1}=\Gamma_{k+1}', S_{k+1}=S_{k+1}'$.


Consider now that $f_{l+k_l+1}$ or its inverse has the label $t$. Put then $\Gamma_{k+1}'=\delta_n(\Gamma_k, R)$ where $R=A(\alpha(f_{l+k_l+1}))$ if $f_{l+k_l+1}$ is positively oriented and $R=B(\alpha(f_{l+k_l+1}))$ otherwise. Put also $Q=A, P=B$ if $R$ is the Cayley graph of $A$ and $Q=B, P=A$ otherwise. We may consider that if $f_{l+k_l+1}$ is positively oriented then $A(\omega(S_k^1))=A(\omega(f_{l+l_k}))$ or $B(\omega(S_k^1))=B(\omega(f_{l+l_k}))$ otherwise. Consider the path $S_{k+1}'=S_k^1\cup f$ where $f$ is the edge going from $\omega(S_k)^1$ and going into the vertex of the second copy of $\Gamma_k$. According to the construction of $\Gamma_{k+1}'$ it is possible to deduce that the length of $S_{k+1}'$ is greater or equal than $k+1$ and the label of $f$ or its inverse equals $t$, the $u$-cycle going from the vertex $\alpha(T)^1$ is maximal in $\Gamma_{k+1}'$ and has the subpath $G$-near to $S_{k+1}'$. Let $Z$ be the representative of an arbitrary $v$-cycle which is maximal in $\Gamma_k$ and goes away from the vertex $p$ and corresponds to the chosen notation of $v$, $Z=z_1...z_q$, the path $z_{j+1}...z_{j+l_k}$ is $G$-near to $S_k$. Then the representative $Z'$ of a $v$-cycle of $\Gamma_{k+1}'$ going from the vertex $p^1$ and corresponding to the chosen notation of $v$ has the subpath $z_{j+1}^1...z_{j+l_k}^1$ which is $G$-near to the path $S_k^1$.

Suppose that $Q(\omega(z_{j+l_k}^1))\neq Q(\omega(S_k^1))=Q(\omega(f_{l+l_k}))$. In this case all vertices of $Z'$ belong to the first copy of $\Gamma_k$. Because otherwise $Z'$ has 1- or 2-near vertices belonging to $Q(\omega(S_k^1))$ or $P(f)$ like in previous cases. So we may suppose that $Q(\omega(z_{j+l_k}^1))=Q(\omega(S_k^1))=Q(\omega(f_{l+l_k}))$ and from the definition of $G$-near paths it is possible to conclude that the edge $z_{j+l_k+1}'^1$ or its inverse has label $t$ where $z_{j+l_k+1}'^1$ is the edge of $Z'$ passing after $e_{j+l_k}^1$. Similarly it is true that if orientations of $z_{j+l_k+1}$ and $f_{l+l_k+1}$ are different then either $Z'$ belongs to the first copy of $\Gamma_k$ or $Z'$ has 1-near vertices in $R^1$. Thus we may assume that $z_{j+l_k+1}$ and $f_{l+l_k+1}$ have equal orientations and therefore $Z'$ is $G$-near to $S_{k+1}'$ and it is possible to regard that $\Gamma_{k+1}=\Gamma_{k+1}', S_{k+1}=S_{k+1}'$.

Notice that if $k=0$ we just consider that $f_{l_k+1}$ is the edge going from $p$ which the first edge of the representative of a $u$-cycle corresponding to the fixed notation of $u$. Besides in this case it is possible to consider that $e_{l+l_k+1}'^1=f_{l+l_k+1}^1$.

The above arguments show that if all maximal in $\Gamma_{k+1}'$ $v$-cycles has the subpath $G$-near to the path $S_{k+1}'$ then the conditions 2), 4), 5) are true for $\Gamma_{k+1}', S_{k+1}'$ and we may put $\Gamma_{k+1}=\Gamma_{k+1}, S_{k+1}=S_{k+1}'$. Otherwise for some $k$ the length of a maximal $u$-cycle is greater than the length of a maximal $v$-cycle and we put $M=k$ and the problem is solved. Suppose that the graph $\Gamma_k$ with properties 1)--5) exists for each natural $k$. Then we put $M=\infty$. In this case in the graph $\Gamma_0$ the representatives of $u$ and $v$-cycles going from one vertex are $G$-near. Since $|\psi(u)|=|\psi(v)|$ and representatives of $u$- and $v$-cycles have no 2-near vertices we deduce that lengths of elements $u$ and $v$ as elements of HNN extension coincide, $m=n, \varepsilon_i=\delta_i, 1\leqslant i\leqslant m+1$. This conditions involve that $v=uh, h\in A\cup B$ besides $h\in A$ if $\varepsilon_{n+1}=1$ and $h\in B$ if $\varepsilon_{n+1}=-1$. Thus we conclude that $v'=u'q, q\in A\cup B, h=\psi_1\times\psi_2(q)$ and for some $j, j\geqslant0, y_k^{-1}hy_k\in A\cup B, x_k^{-1}qx_k\notin A\cup B, x_k=g_0\prod_{i=1}^jt^{\varepsilon_i}g_iz, y_k=g_0'\prod_{i=1}^jt^{\varepsilon_i}g_i'z', z'=\psi_1\times\psi_2(z), z\in\{1, t^{\varepsilon_{j+1}}\}$. This contradicts to the condition about $\psi_1\times\psi_2$ and the set $\Delta$. Theorem 1 is proven.

Note that we used the properties of the base group and of connected subgroups only for to reduce the problem to the case of a finite base group and when we dealt with elements $u$ and $v$ belonging to one coset on one of the connected subgroups or when $u$ or $v$ belongs to a subgroup conjugate to base group. So we assert that the following theorem is proved.

\textbf{Theorem 3.} Consider the HNN extension $H=\<G, t; t^{-1}at=\varphi(a), a\in A<G\>$ of a finite base group $G$ with stable letter $t$ and isomorphic subgroups $A, \varphi(A)$. Fix elements $u, u$ of $H$ such that $u$ is not conjugate to $v^{\pm1}$ and $u, v$. Let $u', v'$ be a cyclically reduced elements conjugate with $u$ and $v$ correspondingly. If $u', v'$ do not belong to one coset on $A$ or $\varphi(A)$ and neither $u$ nor $v$ belongs to $h^{-1}Gh$ for some $h\in H$ then there exists a homomorphism $\psi$ of $G$ onto a finite group such that the orders of $\psi(u), \psi(v)$ are different.

\section{Free products with commutative subgroups.}

We shall use the symbol $G=(A\ast B; [M, N])$ to denote the free product of $A$ and $B$ with commutative subgroups $M<A, N<B$.

Notice that according to the definition of the free product with commutative subgroups it follows that for each homomorphisms $\varphi: A\rightarrow A', \psi: B\rightarrow B'$ there exists the singular homomorphism $\varphi\ast\psi: G\rightarrow G'=(A'\ast B'; [\varphi(M), \psi(N)])$ such that $\varphi\ast\psi|_A=\varphi, \varphi\ast\psi|_B=\psi$ which can defined as follows. Consider the element $g\in G, g=\prod_{i=1}^ng_i, g_i\in A\cup B, 1\leqslant i\leqslant n$. Then $\varphi\ast\psi(g)=\prod_{i=1}^n\epsilon_i(g_i)$ and if $g_i\in A$ then $\epsilon_i=\varphi$ if $g_i\in B$ then $\epsilon_i=\psi$. This algorithm defines the homomorphism because all defining relations of the group $G$ representation has unit images.

\textbf{Definition.} \textsl{The reduced notation of an element $g$ is the presentation of $g$ in generators $\{A\cup B\}$: $g=\prod_{i=1}^ng_i, 1\leqslant i\leqslant n$ satisfying the following properties:}

\textsl{1) there are no entries of the kind: $g_ig_{i+1},
g_i\in A, g_{i+1}\in A$, $g_ig_{i+1}, g_i\in B, g_{i+1}\in B$ in the product;}

\textsl{2) if $n>2$ then there are no entries of the kind:
$g_ig_{i+1}, g_i, g_{i+1}\in M\cup N$ in the product;}

From the properties of free action graphs of the free product $G=(A\ast B; [M, N])$ with commutative subgroups it is clear how the Cayley graph of $G$ with generators $A\cup B$ look like. This involves that the product which is the reduced notation is not equal to unit in $G$ and the number of factors is the same for each such presentation of a fixed element. This number will be called the \textsl{length} of the element $g$. It is evident that each element of $G$ possesses the reduced notation. Also if an element possesses the reduced notation which does not contain elements from the commutative subgroups then this notation is singular.

The element of $G$ will be called cyclically reduced if in its reduced notation the first and the last elements do not belong to one free factor and do not belong to different commutative subgroups in case the length of an element is greater than 2. Notice that two cyclically reduced elements $u$ and $v$ are conjugate if and only if $u$ is a cyclic transposition of $v$.

Consider the free action graph $\Gamma$ of $G=(A\ast B; [M, N])$. Derive the following rule of construction of a new free action graph of $G$. Fix a natural number $n>1$ and a subgraph $R$ which is the Cayley graph of $M$ with the generating set $\{M\}$. Consider $n$ copies of $\Gamma: \Delta_1,..., \Delta_n$. For each $i, 1\leqslant i\leqslant n$, and for each edge $e$ of $\Gamma$ whose label belong to $A\setminus M$ and whose terminal vertex $p$ belongs to $R$ make the following procedure. Let $q$ be the second terminal vertex of $e$. Delete edges $e^i$ and connect vertices $p^i, q^{i+1}$ by the edge $f_i$ whose label equals the label of $e$. If $e$ goes away from $p$ then $f_i$ goes away from $p^i$ and if $e$ goes into $p$ then $f^i$ goes into $p^i$. The obtained graph is the free action graph of $G$. It will be denoted by $\Delta_n(\Gamma, R)$. We can construct the similar graph if we take the subgraph which is the Cayley graph of $N$ with generating set $\{N\}$ in the capacity of $R$.

\textbf{Definition.} \textsl{We shall say that paths $S=e_1...e_n, T=f_1...f_n$ of the free action graph of $(A\ast B; [M, N])$ are near if for each $i, 1\leqslant i\leqslant n$ either $M(\alpha(e_i))=M(\alpha(f_i))$ or $N(\alpha(e_i))=N(\alpha(f_i))$ and either $M(\omega(e_n))=M(\omega(f_n))$ or $N(\omega(e_n))=N(\omega(f_n))$.}

\section{Proof of theorem 2.}

We shall consider that the elements $u$ and $v$ are cyclically reduced. Due to the condition of theorem 2 it is sufficient to consider the following cases:

$$
1) u=ab, v=a'b', a, a'\in M, b, b'\in N, a\neq q^{-1}a'^{\pm1}q,
q\in A
$$
$$
2) u\notin g^{-1}\<M, N\>g, g\in G, v=ab, a\in M, b\in N
$$
$$
3) u\notin g^{-1}\<M, N\>g, v\notin g^{-1}\<M, N\>g, g\in G
$$

In the first case it is sufficient to consider the natural homomorphism of $G$ onto $A$ and to use the 2-order separability of $A$.

Consider the second case. Fix some reduced notation for $u$: $u=g_1...g_n$. According to the conditions of theorem 2 the groups $A$ and $B$ are correspondingly $M$- and $N$-separable and subgroups separable so there exist homomorphisms $\varphi, \psi$ of $A$ and $B$ correspondingly onto finite groups such that $\varphi((\{g_1,..., g_n\}\cap A)\setminus M)=\varphi(\{g_1,..., g_n\}\cap A)\setminus\varphi(M), \varphi(\{g_1,..., g_n\}\cap A)\cap\{1\}$ is empty, $\psi((\{g_1,..., g_n\}\cap B)\setminus N)=\psi(\{g_1,..., g_n\}\cap B)\setminus\psi(N), \psi(\{g_1,..., g_n\}\cap B)\cap\{1\}$ is empty. Thereby there exists the homomorphism $\epsilon$ of $G$ onto free product $G'$ of groups $\varphi(A), \psi(B)$ with commutative subgroups $\varphi(M), \psi(N)$ such that $\epsilon(u)\notin\cup_{g\in G'}g^{-1}\<\varphi(M), \psi(N)\>g$. Thus the element $\epsilon(u)$ has infinite order. The order of $\epsilon(v)$ is finite. According to [8] the group $G'$ is residually finite. So according to the condition on orders of $\epsilon(u), \epsilon(v)$ there exists a homomorphism of $G'$ onto a finite group such that the orders of images of $\epsilon(u)$ and $\epsilon(v)$ are different.

Consider the third case. Fix reduced notations for elements $u$ and $v: u=g_1...g_n, v=h_1...h_m$. By the supposition $n\geqslant2, m\geqslant2$. Put $A_u=\{g_1,..., g_n\}\cap A, B_u=\{g_1,..., g_n\}\cap B, A_v=\{h_1,..., h_m\}\cap A, B_v=\{h_1,..., h_m\}\cap B, A_{uv}=A_u\cup A_v\cup\{a_ia_j^{-1}, a_l'a_k'^{-1}, a_ia_k'^{-1}| a_i, a_j\in A_u, a_k', a_l'\in A_v\}, B_{uv}=B_u\cup B_v\cup\{b_{i'}b_{j'}^{-1}, b_{l'}'b_{k'}'^{-1}, b_{i'}b_{k'}'^{-1}| b_{i'}, b_{j'}\in B_u, b_{k'}', b_{l'}'\in B_v\}$. Since $A$ and $B$ are subgroup separable there exist homomorphisms $\psi_1, \psi_2$ of $A$ and $B$ correspondingly onto finite groups such that $\psi_1(A_{uv}\setminus M)=\psi_1(A_{uv})\setminus\psi_1(M), \psi_1(A_{uv}\setminus\{1\})=\psi_1(A_{uv})\setminus\{1\}, \psi_2(B_{uv}\setminus N)=\psi_2(B_{uv})\setminus\psi_2(N), \psi_2(B_{uv}\setminus\{1\})=\psi_2(B_{uv})\setminus\{1\}$. Thus there exists the homomorphism $\epsilon$ of $G$ onto a group $G'$ which is the free product of groups $A'=\psi_1(A)$ and $B'=\psi_2(B)$ with commutative subgroups $M'=\psi_1(M), N'=\psi_2(N)$ such that the elements $u'=\epsilon(u), v'=\epsilon(v)$ do not belong to subgroups which are conjugate to $A', B', \<M', N'\>$ and $u'$ is not conjugate to $v'^{\pm1}$. In [9] it was proves that the group $G'$ is cyclic subgroup separable hence there exists the homomorphism $\varphi_1$ of $G'$ onto a finite group such that elements not belonging to $\<u'\>$ and $\<v'\>$ which are equal to $(g_ig_{i+1}...g_n)^{-1}qg_1...g_j, (h_kh_{k+1}...h_m)^{-1}qh_1...h_l, q=a'b', a'\in A', b'\in B',$ have images after $\varphi_1$ that do not belong $\<\varphi_1(u)\>, \<\varphi_1(v)\>$. Particularly all words of length 1 or 2 are not mapped into unit by $\varphi_1$. So in the Cayley graph $\Lambda$ of $G'$ with the generating set $\{A'\cup B'\}$ the representatives of $u'$- and $v'$-cycles have no 2-near vertices. Let $g_i'=\epsilon(g_i), h_j=\epsilon(h_j), 1\leqslant i\leqslant\leqslant n, 1\leqslant j\leqslant m$.


Suppose that $\varphi_1(u')$ and $\varphi_1(v')$ have equal orders. Consider the inductive process of construction of finite free action graphs $\Gamma_n$ of the group $G'$ satisfying the following properties ($0\leqslant n\leqslant r, r$ will be chosen later):

1) the length of each $u'$-cycle divides the length of the maximal $u'$-cycle; the same is true for $v'$-cycles;

2) representatives of $u'$- and $v'$-cycles have no 2-near vertices

3) lengths of maximal $u'$- and $v'$-cycles coincide

4) in $\Gamma_n$ when $n>0$ there exists the path $S_n$ whose length is greater or equal than $n$ such that there exists the representative of a maximal $u'$-cycle having the subpath near to $S_n$ and each representative of each maximal $v'$-cycle has a subpath near to $S_n$;

5) the graph $\Gamma_n$ has no a cycle $e_1e_2e_3e_4$ such that Lab$(e_1)$, Lab$(e_3)\in A'$, Lab$(e_2)$, Lab$(e_4)\in B'$, and at least one of the elements from Lab$(e_1)$, Lab$(e_3)$ does not belong to $M'\times N'$;

According to the definition of the Cayley graph the graph $Cay(\varphi_1(G'); A'\cup B')$ may be regarded as $\Gamma_0$.

Construct the graph $\Gamma_1$. Let $m_1$ be the length of a $u$-cycle in $\Gamma_0$. Fix an arbitrary vertex $p$ in $\Gamma_0$. We may assume that first letters in the chosen notations of $u$ and $v'$ do not belong to $M'\times N'$. Consider the representative $T=e_1...e_s$ of a $u$-cycle starting at $p$ and corresponding to the chosen notation of $u$. Let $e_{1+w}$ be the first edge of $T$ following after $e_1$ whose label does not belong to $M'\times N'$. According to the definition of the reduced notation $w\in\{1, 2\}$. Suppose that Lab$(e_1)\in A'\setminus M'$. Put $\Gamma_1'=\delta_{m_1^2}(\delta_{m_1}(\Gamma_0; M'(\omega(e_1))); R)$ where $R=N'(\omega(e_2^2))$ if $w=1$ and $R=M'(\omega(e_3^2))$ if $w=2$. Put $\Delta'=\delta_{m_1}(\Gamma_0; M'(\omega(e_1)))$ and $e'$ is the edge of the first copy of $\Delta'$ which is appended instead of the edge $e_1^1$, $f'$ is the edge of the first copy of the graph $\Delta'$ which is appended instead of $e_{1+w}^2$. The graph $\Delta'$ satisfies the property 2) (because otherwise the property is not held for $\Gamma_0$). Since there are no 1- and 2-near vertices in representatives of $u'$- and $v'$-cycles of $\Gamma_0$ we may deduce that either each representative of each $v'$-cycle in $\Delta'$ belongs just to one copy of $\Gamma_0$ or the representative of each maximal in $\Delta'$ $v'$-cycle passes through an edge with label from $A'\setminus M'$ going into the vertex of $M'(\alpha(e_2^2))$. If the latter is not true then $r=0$ and in the graph $\Delta'$ the length of a maximal $u'$-cycle exceeds the length of a maximal $v'$-cycle (we use here the property 3) for $\Gamma_0$ and 1) for $\Delta'$). Otherwise we have two alternatives for $\Gamma_1'$. Let $r$ be the vertex of $\Delta'$ which is the beginning of the maximal in $\Delta'$ $v'$-cycle. Consider its representative $R=y_1...y_s$. Suppose that the edge $y_i$ has the label from $A'\setminus M'$ and goes into the vertex from $M'(\omega(e'))$. Suppose that $w=1$. Then if $\alpha(y_{i+1})\neq \alpha(f')$ we have that all vertices of representatives of the $v'$-cycle $R'$ of $\Gamma_1'$ going from $r_1$ belong to the first copy of $\Delta'$ because otherwise $R$ has 2-near vertices from $A'(\alpha(e')\cup B'(\alpha(f'))$. In the other case representatives of $R'$ pass through the edge $f'$. So it is proven that if there are no such $v'$-cycles then in $\Gamma_1'$ the length of a maximal $u'$-cycle is greater than the length of a maximal $v'$-cycle. Otherwise we have the graph $\Gamma_1=\Gamma_1'$ with properties 1)--5), $S_1=f'$. Consider that $w=2$. If the edge $y_{i+w}$ does not go into the vertex of $M'(\alpha(f'))$ then all the vertices of representatives of $R'$ belong to one copy of $\Delta'$ since otherwise representatives of $R'$ would have 2-near vertices in $A'(\alpha(f'))$. In this case as earlier the length of a maximal $v'$-cycle is less than the length of a maximal $u'$-cycle in $\Gamma_1'$ and theorem is proven. Otherwise $y_{i+w}$ is near to $f'$ (due to the property 5) for $\Gamma_1'$) and there exists the graph $\Gamma_1=\Gamma_1'$ with properties 1)--5) and $S_1=f'$.

Consider there exists the graph $\Gamma_k, k>1$. Construct the graph $\Gamma_{k+1}$. Let $n_1$ be the length of a maximal $u'$-cycle in $\Gamma_k$. Consider the representative of such $u'$-cycle $T=e_1...e_d$ such that the path $e_{i+1}...e_{i+k'}$ is near to the path $S_k=f_1...f_{k'}$. Consider also an arbitrary representative of a maximal $v'$-cycle $R=r_1...r_s$. The path $r_{l+1}...r_{l+k'}$ is near to the path $S_k$. Let $e_{i+k'+w}$ be the first edge of $T$ following after $e_{i+k'}$ whose label does not belong to $M'\times N'$. Notice that according to the definition of the reduced notation $w\in\{1, 2\}$. We may assume without loss of generality that Lab$(e_{i+k'+w})\in A'\setminus M'$. Construct the graph $\Gamma_{k+1}'=\Delta_{n_1}(\Gamma_k; M'(\omega(e_{i+k'+w})))$. The graph $\Gamma_{k+1}'$ satisfies the condition 1) because of the condition on $n_1$. Conditions 2) and 5) are held because otherwise they are not held for $\Gamma_k$.

Suppose that the edge $r_{l+k'+w}$ does not go into the vertex from $M'(e_{i+k'+1})$. In this case all vertices of a $v'$-cycle going from the vertex $\alpha(r_1)^1$ belong to the first copy of the graph $\Gamma_k$ because otherwise representatives of this $v'$-cycle have 1-near vertices if $w=1$ or 2-near vertices in case $w=2$. If this condition is true for some maximal in $\Gamma_k$ $v'$-cycles the length of a maximal in $\Gamma_{k+1}'$ $v'$-cycle is less than the length of a maximal $u'$-cycle. In this case we put $r=k$ and theorem is proven. So we assume that the representative $R$ in $\Gamma_k$  goes into vertices of $M'(\omega(e_{i+k'+w}))$ and hence due to the property 4) for $\Gamma_k$ and 5) for $\Gamma_{k+1}'$ we deduce that the path $e_{i+1}...e_{i+k'+w}$ is near to the path $r_{l+1}...r_{l+k'+w}$. So there exists the extension $S_{k+1}$ in $\Gamma_{k+1}'$ of $S_k^1$ and the property 4) is held. That is in the graph $\Gamma_{k+1}'$ each representative of each maximal $v'$-cycle has the subpath near to the path $S_k^1$ and contains the edge with label from $A'\setminus M'$ going into the vertex of $M'(\omega(e_{i+k'+w}))^2$. Particularly $\Gamma_{k+1}'$ satisfies the condition 3). Also $\Gamma_{k+1}'$ satisfies the condition 4). If $w=1$ the $S_k^1\cup e'$ is taken in the capacity of $S_{k+1}$ where $e'$ is an arbitrary edge with label from $A'\setminus M'$ going from the vertex $\omega(S_k^1)$ into the vertex $\alpha(e_{i+k'+w+1}^2)$. If $w=2$ we put $S_{k+1}=S_k^1\cup e_1'\cup e_2'$ where $e_1'$ is an edge with label from $N'$ going from $\omega(S_k^1)$ into the vertex of the subgraph $M'(\alpha(e_{i+k'+w}^1)), e_2'$ is an edge with label of $A'\setminus M'$ going from the vertex $\omega(e_1')$ and going into the vertex $\omega(e_{i+k'+w}^2)$. So we obtain the graph $\Gamma_{k+1}$.

Consider that the graph $\Gamma_k$ exists for each natural number $k$. It means that in the graph $\Gamma_0$ the representatives of $u'$- and $v'$-cycles going from one vertex are near. Since representatives of $u'$- and $v'$-cycles have no 2-near vertices in $\Gamma_0$ then lengths of $u'$ and $v'$ are equal and $v'=u'^{\pm1}$. Violation. This involves that the representatives of $u'$- and $v'$-cycles can not be near in $\Gamma_0$. Theorem 2 is proven.


\begin{center}
\large{Acknowledgements.}
\end{center}

The author thanks A. A. Klyachko for valuable comments and discussions.

\begin{center}
\large{References.}
\end{center}

1. \textsl{Klyachko A. A.} Equations over groups, quasivarieties,
and a residual property of a free group // J. Group Theory. 1999.
\textbf{2}. 319--327.

2. \textsl{Wise, Daniel T.} Subgroup separability of graphs of free groups with cyclic edge groups.
Q. J. Math. 51, No.1, 107-129 (2000). [ISSN 0033-5606; ISSN 1464-3847]

3. \textsl{A. Minasyan, D. Osin}. Normal automorphisms of relatively hyperbolic groups. To appear, arXiv:0809.2408

4. \textsl{Graham A. Niblo.} H.N.N. extensions of a free group by \textbf{z} which are subgroup separable. Proc. London Math. Soc. (3), 61(1):18--23, 1990.

5. \textsl{Yedynak V. V.} Separability with respect to order. \textsl{Vestnik Mosk. Univ. Ser. I Mat. Mekh.} \textbf{3}: 56-58, 2006.

6. \textsl{Yedynak V. V.} Order separability of free product of free groups with cyclic amalgamation. arXiv:1007.0157. This work was submitted in \textsl{Vestnik Mosk. Univ. Ser. I Mat. Mekh.}

7. \textsl{Lyndon, R. C., Schupp, P. E.} (1977). Combinatorial group
theory.\ \textsl{Springer-Verlag}.

8. \textsl{Loginova Ye. D.} Residual finiteness of the free product of two groups with commutative subgroups // Sibir. Math. J. 1999. Ò. 40, \textbf{2}, 395-407

9. \textsl{Loginova Ye. D.} The cyclic subgroup separability of the free product of two groups with commutative subgroups // Vestnik Ivanovsk. Gos. Univ, \textbf{3}, (2000), 47-53

\end{document}